\documentclass[reqno,a4paper,11pt]{amsart}
\usepackage{hyperref}
\usepackage{amsthm,amssymb}
\usepackage{graphicx}
\usepackage{enumerate}

\theoremstyle{definition}
\newtheorem{definition}{Definition}

\theoremstyle{remark}

\theoremstyle{plain}

\newtheorem{proposition}[definition]{Proposition}
\newtheorem{theorem}[definition]{Theorem}
\newtheorem{cor}[definition]{Corollary}

\newcommand{\set}[1]{\left\{{#1}\right\}}

\newcommand{\Z}{{\mathbb{Z}}}

\newcommand{\CompOf}[2]{#1 \vDash #2}

\newcommand{\w}[1]{{\lvert #1 \rvert}}
\newcommand{\NC}{\mathfrak{N}}

\begin{document}
\title{Standard paths in another composition poset}
\author{Jan Snellman}
\address{Department of Mathematics\\
Stockholm University\\
SE-10691 Stockholm, Sweden}
\email{Jan.Snellman@math.su.se}
\begin{abstract}
  Bergeron, Bousquet-M{\'e}lou and Dulucq \cite{StPa} enumerated  paths
  in the Hasse diagram of the following poset: the underlying set is
  that of all compositions, and a composition \(\mu\) covers another
  composition \(\lambda\) if \(\mu\) can be obtained from \(\lambda\)
  by adding \(1\) to one of the parts of \(\lambda\), or by inserting
  a part of size \(1\) into \(\lambda\).

We employ the methods they developed in order to study the same
problem for the following poset: the underlying set is the same, but
\(\mu\) covers \(\lambda\) if \(\mu\) can be obtained from \(\lambda\)
  by adding \(1\) to one of the parts of \(\lambda\), or by inserting
  a part of size \(1\) at the left or at the right of
  \(\lambda\). This poset is of interest because of its relation to
  non-commutative term orders \cite{Snellman:Ncterm}.
\end{abstract}

\maketitle

\begin{section}{Definition of standard paths}
  By a \emph{composition} \(P\) we mean a sequence of positive
  integers \((p_1,p_2,\dots,p_k)\), which are the \emph{parts} of
  \(P\). We define the \emph{length} \(\ell(P)\) of \(P\) as the
  number of parts, and the \emph{weight} \(\w{P} = \sum_{i=1}^k p_k\)
  as the sum of its parts. If \(P\) has weight \(n\) then \(P\) is a
  composition of \(n\), and we write \(\CompOf{P}{n}\).

  We say that a composition \(Q\) covers a composition \(P\) if \(Q\)
  is obtained from \(P\) either by adding 1 to a part of \(P\), or by
  inserting a part of size 1 to the left, or by inserting a part of
  size 1 to the right. Thus, \(P=(p_1,p_2,\dots,p_k)\) is covered by 
  \begin{enumerate}
  \item \((1,p_1,\dots,p_k)\),
  \item \((p_1,\dots,p_k,1)\),
  \item and, for \(1 \le i \le k\), \((p_1,\dots, p_i +1, \dots,p_k)\).
  \end{enumerate}

  Extending this relation by transitivity makes the set of all
  compositions into a partially  ordered 
  set, which we denote by \(\NC\). This is in accordance with the
  notations in the author's article \textbf{A poset classifying
    non-commutative term orders} \cite{Snellman:Ncterm}, where \(\NC\)
  was used for the 
  following   isomorphic poset of words: the underlying set is
  \(X^*\), the free associative monoid on
  \(X=\set{x_1,x_2,x_3,\dots}\), and \(m_1=x_{i_1} \cdots x_{i_r}\) is
  smaller than \(m_2\) if \(m_2\) can be obtained from \(m_1\) by
  a sequence of operations of the form 
  \begin{enumerate}[(i)]
  \item Multiply by a word to the left,
  \item Multiply by a word to the right,
  \item Replace an occuring \(x_i\) with an \(x_j\), with \(j > i\).
  \end{enumerate}
  The bijection \((p_1,p_2,\dots,p_k) \mapsto x_{p_1} \cdots x_{p_k}\)
  is an order isomorphism between these two partially ordered sets.

  On the other hand, the partial order \(\Gamma\) on compositions
  studied by Bergeron, Bousquet-M{\'e}lou and Dulucq in
  \textbf{Standard paths in the composition poset}
  \cite{StPa} is different, since in \(\Gamma\) the composition
  \(P=(p_1,p_2,\dots,p_k)\) is covered by 
\begin{enumerate}
  \item \((1,p_1,\dots,p_k)\),
  \item \((p_1,\dots,p_k,1)\),
  \item for \(1 \le i \le k\), \((p_1,\dots, p_i +1, \dots,p_k)\),
  \item for \(1 \le i < k\), \((p_1,\dots, p_i ,1,p_{i+1}, \dots,p_k)\).
  \end{enumerate}
  
  \(\Gamma\) and \(\NC\) coincide for compositions of weight \(\le
  4\). In Figure~\ref{fig:HasseNC4} this part of the Hasse diagram is
  depicted. We have that \((2,2) \le (2,1,2)\) in \(\Gamma\) but not
  in \(\NC\), so the rest of the respective Hasse diagrams differ.

  \setlength{\unitlength}{1cm}
  \begin{figure}[htbp]
    \begin{center}
      \begin{picture}(12,6) 
        \put(5.9,0.7){\(0\)}
        \put(6,1){\line(0,1){1}}
        
        \put(5.3,2){\(1\)}
        \put(6,2){\line(-2,1){2}}
        \put(6,2){\line(2,1){2}}

        \put(3.1,2.6){\(11\)}
        \put(4,3){\line(-3,1){3}}
        \put(4,3){\line(1,1){1}}
        \put(4,3){\line(3,1){3}}

        \put(8.5,2.7){\(2\)}
        \put(8,3){\line(2,1){2}}
        \put(8,3){\line(-1,1){1}}
        \put(8,3){\line(-3,1){3}}

        \put(0.3, 4){\(111\)}
        \put(1,4){\line(-1,1){1}}
        \put(1,4){\line(1,1){1}}
        \put(1,4){\line(3,1){3}}
        \put(1,4){\line(5,1){5}}

        \put(4,4){\(12\)}
        \put(5,4){\line(-1,1){1}}
        \put(5,4){\line(1,1){1}}
        \put(5,4){\line(2,1){2}}
        \put(5,4){\line(3,1){3}}

        \put(6,4){\(21\)}
        \put(7,4){\line(-5,1){5}}
        \put(7,4){\line(-3,1){3}}
        \put(7,4){\line(0,1){1}}
        \put(7,4){\line(3,1){3}}

        \put(10.5,4){\(3\)}
        \put(10,4){\line(-2,1){2}}
        \put(10,4){\line(0,1){1}}
        \put(10,4){\line(2,1){2}}

        \put(0,5.3){\(1111\)}

        \put(2,5.3){\(211\)}

        \put(4,5.3){\(121\)}

        \put(6,5.3){\(112\)}
        
        \put(7,5.3){\(22\)}
        
        \put(8,5.3){\(13\)}

        \put(10,5.3){\(31\)}

        \put(12,5.3){\(4\)}
      \end{picture}
      \caption{The Hasse diagram of \(\NC\).}
      \label{fig:HasseNC4}
    \end{center}
  \end{figure}

  Following \cite{StPa} we define a \emph{standard path of length
  \(n\)} to be a sequence \(\gamma=(P_0,P_1,P_2,\dots,P_n)\) of
  compositions such that 
  \begin{equation}
    \label{eq:spath}
    P_0 \prec P_1 \prec P_2 \prec \cdots \prec P_n, \qquad \CompOf{P_i}{i}.    
  \end{equation}
  The partial order is now that of \(\NC\). For instance,
  \begin{equation} \label{eq:sp} 
  \rho=((),(1),(1,1),(1,2),(1,1,2))
  \end{equation}
   is a standard path of length
  4, corresponding to a saturated chain in Hasse diagram of \(\NC\)
  between the minimal element 
  \(()\) and the element \((1,1,2)\).

  We furthermore define the \emph{diagram} of a composition
  \(P=(p_1,\dots,p_k)\) to be the set of points \((i,j) \in \Z^2\)
  with \(1 \le j \le p_i\). Alternatively, we can replace the node
  \((i,j)\) by the square with corners
  \((i-1,j-1)\),\((i-1,j)\),\((i,j-1)\) and \((i,j)\). So the
  composition \((1,1,2)\) has diagram 
    \setlength{\unitlength}{0.4cm}
    \begin{picture}(3.5,2.5)
    \put(0,0){\line(1,0){3}}
    \put(0,1){\line(1,0){3}}
    \put(2,2){\line(1,0){1}}
    \put(0,0){\line(0,1){1}}
    \put(1,0){\line(0,1){1}}
    \put(2,0){\line(0,1){2}}
    \put(3,0){\line(0,1){2}}
  \end{picture}.
  For a standard path \(\gamma=(P_1,\dots,P_n)\) ending at \(P_n\) we
  label the boxes in the diagram of \(P_n\) in the order that they
  appear in the path. To avoid ambiguity, we use the convention that
  whenever \(P_i\) consists of \(i\) ones and \(P_{i+1}\) consists of
  \(i+1\) ones, the extra one is considered to have been added to the
  left. So for the path \(\rho\)  the corresponding
  tableau is 
    \begin{picture}(3.5,2.5)
    \put(0,0){\line(1,0){3}}
    \put(0,1){\line(1,0){3}}
    \put(2,2){\line(1,0){1}}
    \put(0,0){\line(0,1){1}}
    \put(1,0){\line(0,1){1}}
    \put(2,0){\line(0,1){2}}
    \put(3,0){\line(0,1){2}}
    \put(0.2,0.2){4}
    \put(1.2,0.2){2}
    \put(2.2,0.2){1}
    \put(2.2,1.2){3}
  \end{picture}.

  Clearly, two different standard paths give rise to different
  tableau. Furthermore, the tableau that occurs as tableau of standard
  paths must be increasing in every column, and have the additional
  property that 
  whenever the numbers \(1,2,\dots,k\) occur as a contiguous sequence
  on the bottom row, then that sequence is \(k,k-1,\dots,2,1\). This
  is a necessary but not sufficient condition.

  The underlying diagram of a tableau is called its \emph{shape}, and
  we define the shape of a standard path to be the shape of its
  tableau. We define the \emph{height} and \emph{width} of a diagram
  to be the height and width of the smallest rectangle containing
  it. Hence, the standard path \(\rho\)  has width 2 and height
  1.

\end{section}

  \begin{section}{Enumeration of standard paths of fixed width}
    Let \(\NC_{(k)}\) denote the subposet of compositions of width
    \(k\). For a path \(\gamma\) of shape \((p_1,p_2,\dots,p_k)\) we
    set 
    \begin{equation}
      \label{eq:vwidth}
      v(\gamma)=x_1^{p_1} x_2^{p_2} \cdots x_k^{p_k}
    \end{equation}
    We want to compute the generating function
    \begin{equation}
      \label{eq:fk}
      f_k(x_1,\dots,x_k) = \sum_{\gamma \text{ path of width } k}
      v(\gamma) 
    \end{equation}
    
    \begin{theorem}\label{thm:fk}
      The generating funktion \(f_k(x_1,\dots,x_k)\) of standard paths
      of widht \(k\) is a rational function given by the following
      recursive relation:
      \(f_0 =1\), \(f_1(x_1) = x_1 (1-x_1)^{-1}\), and for \(k > 1\)
      \begin{equation}
        \label{eq:fkrel}
        f_k(x_1,\dots,f_k) = \frac{
          x_1f_{k-1}(x_2,\dots,x_k) +
          x_kf_{k-1}(x_1,\dots,x_{k-1}) - x_1 \cdots x_k 
        }
        {1-x_1 - \ldots -x_k} 
      \end{equation}
    \end{theorem}
    \begin{proof}
      A tableau of width \(k\) can be obtained by adding a new cell
      either 
      \begin{enumerate}[-]
      \item at the top of a column of another tableau of width \(k\),
      \item at the beginning  of a tableau of width \(k-1\),
      \item or at the end of a tableau of width \(k-1\).
      \end{enumerate}
      These three cases correspond respectively to
      \begin{math}
        (x_1 +x_2 + \ldots + x_k) f_k
      \end{math},
      to
      \begin{math}
        x_1f_{k-1}(x_2,\dots,x_k)
      \end{math},
      and to
      \begin{math}
        x_kf_{k-1}(x_1,\dots,x_{k-1}).
      \end{math}
      However, if the tableau has shape \((1,\dots,1)\) then the last
      two operations give the same result.
      Hence 
      \begin{displaymath}
        f_k = (x_1 +x_2 + \ldots + x_k) f_k +
        x_1f_{k-1}(x_2,\dots,x_k) + 
        x_kf_{k-1}(x_1,\dots,x_{k-1})
        - x_1 \cdots x_k,
      \end{displaymath}
      from which \eqref{eq:fkrel} follows.
    \end{proof}
    We obtain successively
    \begin{equation}
      \label{eq:fksucc}
      \begin{split}
        f_0 &= 1\\
        f_1 &= \frac{x_1}{1-x_1} \\
        f_2 &= \frac{x_1x_2(1-x_1x_2}{(1-x_1)(1-x_2)(1-x_1-x_2)}\\
        f_3 &= \Bigl({x_{{1}}}^{2}{x_{{2}}}^{2}x_{{3}} + {x_{{1}}}^{2}
x_{{2}}{x_{{3}}}^{2} + x_{{1}}{x_{{2}}}^{3}x_{{3}} + x_{{1}}{x_{{2}}}^{2}
{x_{{3}}}^{2} - {x_{{1}}}^{2}{x_{{2}}}^{2} \\
& -
4\,{x_{{1}}}^{2}x_{{2}}x_{{3}} - {x_{{1}}}^{2}{x_{{3}}}^{2} -
x_{{1}}{x_{{2}}}^{3} - 7\,x_{{1}}{x_{{2}}}^{2}x_{{3}} -
4\,x_{{1}}x_{{2}}{x_{{3}}}^{2} - \\
& {x_{{2}}}^{3}x_{{3}} -
{x_{{2}}}^{2}{x_{{3}}}^{2} 
+ 2\,{x_{{1}}}^{2}x_{{2}} + 2\,{x_{{1}}}^{2}x_{{3}} 
+
5\,x_{{1}}{x_{{2}}}^{2} + 12\,x_{{1}}x_{{2}}x_{{3}} \\
& +
2\,x_{{1}}{x_{{3}}}^{2} + {x_{{2}}}^{3} + 5\,{x_{{2}}}^{2}x_{{3}} +
2\,x_{{2}}{x_{{3}}}^{2} - 5\,x_{{1}}x_{{2}} - 4\,x_{{1}}x_{{3}} \\
& -3\,{x_{{2}}}^{2} - 5\,x_{{2}}x_{{3}} 
+ x_{{2}} + 1 
\Bigr ) \\
&\times x_{{1}}x_{{2}}x_{{3}} \\
&\times \left (1-x_{{1}}\right )^{-1} \left (1-x_{{2}}\right )^{-1}\left
  (1-x_{{3}}\right )^{-1}
\left (1-x_{{1}}-x_{{2}}\right)^{-1} 
\left (1-x_{{2}}-x_{{3}}\right)^{-1} \\
& \times
\left (1 - x_{{1}} - x_{{2}} - x_{{3}}\right )^{-1}
      \end{split}
    \end{equation}

\begin{theorem}\label{thm:fk}
  For each \(k\),
  \begin{equation}
  \label{eq:fkstructure}
  f_k(x_1,\dots,x_k) = \frac{x_1\cdots x_k }{
    \prod_{i=1}^k \prod_{j=i}^k (1-x_i -x_{i+1} - \ldots -x_{j})
    } \tilde{f}_k(x_1,\dots,x_k)
\end{equation}
where \(\tilde{f}_k\) is a polynomial.  
\end{theorem}
\begin{proof}
  This is true for \(k=0,1\). Assume that \(f_{k-1}\) has the above form.
  Then 
  \begin{multline}
      f_k (1-x_1 - \cdots -x_k) =  
          x_1f_{k-1}(x_2,\dots,x_k) +
          x_kf_{k-1}(x_1,\dots,x_{k-1}) - x_1 \cdots x_k \\
          = x_1 x_2 \cdots x_k \tilde{f}_{k-1}(x_2,\dots,x_k) 
          \prod_{i=2}^k \prod_{j=i}^k (1-x_i-\cdots -x_j)^{-1}
          \\
           +x_k x_1 \cdots x_{k-1} \tilde{f}_{k-1}(x_1,\dots,x_{k-1}) 
          \prod_{i=1}^{k-1} \prod_{j=i}^{k-1} (1-x_i-\cdots -x_j)^{-1}
          - x_1\cdots x_k
  \end{multline}
  hence 
  \begin{multline}
    \frac{f_k (1-x_1 - \cdots -x_k) \prod_{i=1}^k \prod_{j=i}^k (1-x_i
      - \cdots -x_j)}{x_1 \cdots x_k} \\
    = \tilde{f}_{k-1}(x_2,\dots,x_k)
  \prod_{j=1}^k (1-x_1 - \cdots -x_j)
+  \tilde{f}_{k-1}(x_1,\dots,x_{k-1}) \prod_{i=1}^k (1-x_i - \cdots
-x_k)
\\- \prod_{i=1}^k \prod_{j=i}^k(1-x_i - \cdots - x_j)
  \end{multline}
\end{proof}

Let \(a_{n,k}\) denote the number of standard paths of width \(k\) and
length \(n\), and let 
\begin{displaymath}
  L_k(t) = \sum_{n\ge 0} a_{n,k}t^n
\end{displaymath}
be the generating function for the number of
standard paths of width \(k\) and lenght \(n\).
Then \(L_k(t) = f_k(t,\dots,t)\).  This substitution results in some
cancellation in the numerator and denominator; we have that 
\begin{equation}
  \label{eq:Lk}
  \begin{split}
    L_1(t) &= \frac{t}{1-t} \\
    L_2(t) &= \frac{t^2(1+t)}{(1-t)(1-2t)}\\
    L_3(t) &= \frac{t^3(1+5t-2t^2)}{(1-t)(1-2t)(1-3t)} \\
    L_4(t) &= \frac{t^4(1+16t-15t^2 + 6t^3)}{(1-t)(1-2t)(1-3t)(1-4t)}
    \\
    L_5(t) &= \frac{t^5(1+42t-65t^2 + 62t^3
      -24t^4)}{(1-t)(1-2t)(1-3t)(1-4t)(1-5t)} \\ 
  \end{split}
\end{equation}

\begin{proposition}
\begin{equation}
  \label{eq:Lks}
  L_k(t) = \frac{t^k\tilde{L}_k(t)}{\prod_{i=1}^k(1-it)}
\end{equation}
where \(\tilde{L}_k(t)\) is a polynomial of degree \(k-1\) with 
\(\tilde{L}_k(1)=2^{k-1}\).
\end{proposition}
\begin{proof}
  The recursive relation \eqref{eq:fkrel} specializes to 
  \begin{equation}
    \label{eq:fkrel2}
    L_k = \frac{2tL_{k-1} -t^k}{1-kt} 
  \end{equation}
  Assume \eqref{eq:Lks} for a fixed \(k\); then \eqref{eq:fkrel2}
  gives
  \begin{displaymath}
    L_{k+1} = \frac{2t^{k+1} \tilde{L}_k -t^{k+1} \prod_{i=1}^k(1-it)}
    {\prod_{i=1}^{k+1}(1-it)}
  \end{displaymath}
  Since \(\tilde{L_k}\) has degree \(k-1\) and evaluates to
  \(2^{k-1}\) at 1, we get that \(\tilde{L}_{k+1} =  2 \tilde{L}_k -
  \prod_{i=1}^k(1-it)\) has degree \(k\) and evaluates to \(2^k\) at
  1. The assertion now follows by induction.
\end{proof}

\begin{cor}
  For a fixed \(k\), 
  \begin{equation}
    \label{eq:ank}
    a_{n+k,k} \sim \frac{k^{k-1}}{(k-1)!} k^n \qquad \text{ as } n \to \infty
  \end{equation}
\end{cor}
\begin{proof}
  This follows from the previous Proposition, and the partial
  fraction decomposition 
  \begin{equation}
    \label{eq:pf}
    \begin{split}
    \prod_{i=1}^k (1-it)^{-1} &= \sum_{j=1}^k v_{j,k} (1-it)^{-1} \\
    v_{k,k} &= \frac{k^{k-1}}{(k-1)!}
    \end{split}
  \end{equation}
\end{proof}
\end{section}

  \begin{section}{Enumeration of standard paths of bounded height}
    Let \(\NC_{n,i,j}^{(k)}\) denote the poset of compositions of
    \(n\) with height \(\le k\), having \(i\) parts of size 1 and
    \(j\) parts of size \(\ge 2\). Let \(\gamma_{n,i,j}^{(k)}\) be the
    number of standard paths with endpoint in \(\NC_{n,i,j}^{(k)}\).

    We will derive a recurrence relation for \(\gamma_{n,i,j}^{(2)}\).
    Note that a tableau of height \(\le 2\), with \(i\) parts of
    size \(1\) and \(j\) parts of size \(2\), 
    has a total of \(n=i + 2j\) boxes, so \(\gamma_{n,i,j}^{(2)}=0\)
    unless \(n=i+2j\). 
    Put \(c_{i,j}^{(2)} = \gamma_{i+2j,i,j}^{(2)}\). A tableau  with
    \(i\) parts of 
    size \(1\) and \(j\) parts of size \(2\),
    can be obtained
    \begin{enumerate}[-]
    \item from a tableau with \(i-1\) parts of size 1 and \(j\) parts
      of size \(2\), by adding a part of size 1 to the left,
    \item or from a tableau with \(i-1\) parts of size 1 and \(j\) parts
      of size \(2\), by adding a part of size 1 to the right,
    \item or from a tableau with \(i+1\) parts of size 1 and \(j-1\) parts
      of size \(2\), by adding a box to a part of size 1.
    \end{enumerate}
    For the composition consisting of \(n\) ones the first
    two ways are identical, which gives the recurrence
    \begin{equation}
      \label{eq:hrec}
      \begin{split}
      \gamma_{n,i,j}^{(2)} &= 
      2\gamma_{n-1,i-1,j}^{(2)} + (i+1)\gamma_{n-1,i+1,j-1}^{(2)}
      -\delta_j^0 \\
      c_{i,j}^{(2)} &= 
      2c_{i-1,j}^{(2)} + (i+1)c_{i+1,j-1}^{(2)} -\delta_j^0 
      \end{split}
    \end{equation}
    where \(\delta_i^j\) is the Kronecker delta.
    We get that \(c_{n,0}^{(2)} =\gamma_{n,n,0}^{(2)}= 1\),
    \(\gamma_{n,i,0}^{(2)}= 0\) for \(i \neq n\). For small values of
    \(i,j\), \(c_{i,j}^{(2)}\) is as in table~\ref{tab:c2}

    \begin{table}[htbp]
      \centering
    \begin{tabular}{|rr|rrrrrrrr|}
      \hline 
       & j & 0 & 1 & 2 & 3 & 4 & 5 & 6 & 7     \\

      i&   &   &   &   &   &   &   &   &       \\ \hline

      0&   &  1&  1&  4& 30& 336& 5040& 95040& 2162160\\

      1&   &1 &  4 &  30 &  336 &  5040 &  95040 &  2162160 &
      57657600 \\ 

      2&   & 1& 11& 138& 2184& 42480& 986040& 26666640& 824503680\\

3&  &1 & 26 & 504 & 10800 & 265320 & 7447440 & 236396160 & 8393898240 \\ 

4&  &1 & 57 & 1608 & 45090 & 1368840 & 45765720 & - & -\\ 

5&  &1 & 120 & 4698 & 167640 & 6174168 & 242686080 & - &
- \\ 

6&  &1 & 247 & 12910 & 572748 & 25192440 & 1151011680 & - &
- \\ 

7&  &1 & 502 & 33924 & 1834872 & 95091360 & 4999942080 & - &
- \\ 

8&  &1 & 1013 & 86172 & 5588310 & 337239840 & - & - &
-  \\ 
 \hline
    \end{tabular}
      \caption{Values of \(c_{i,j}^{(2)}\) for small \(i,j\)}
      \label{tab:c2}
    \end{table}

    \begin{theorem}\label{thm:c21}
      Put 
      \begin{equation}
        \label{eq:Pk}
        P_k(x) = \sum_{n =0}^\infty c_{n,k}^{(2)}x^n
      \end{equation}
      Then \(P_0(x) = (1-x)^{-1}\) and 
      \begin{equation}
        \label{eq:Pkrel}
        P_k(x) = \frac{\frac{d}{dx}P_{k-1}(x)}{1-2x}
      \end{equation}
    \end{theorem}
    \begin{proof}
      Since \(c_{n,0}^{(2)}=1\) it follows that 
      \(P_0(x) = \sum_{n=0}^\infty c_{n,0}^{(2)}x^n =(1-x)^{-1}\).
      
      Now, multiply \eqref{eq:hrec} with \(x^i\) and sum over all \(i
      \ge 0\) to 
get that  
      \begin{equation}
        \label{eq:hrec2}
        \sum_{i \ge 0} c_{i,j}^{(2)}x^i = 
        2\sum_{i \ge 1} c_{i-1,j}^{(2)}x^i + \sum_{i \ge
          0}(i+1)c_{i+1,j-1}^{(2)}x^i 
      \end{equation}
which means that 
\begin{equation}
  \label{eq:hrec3}
  P_j(x) = 2xP_j(x) + P'_{j-1}(x) 
\end{equation}
    \end{proof}
    We get that 
    \begin{equation}
      \label{eq:pkval}
      \begin{split}
        P_1(x) & = (1-x)^{-2}(1-2x)^{-1}\\
        P_2(x) & = 2!(1-x)^{-3}(1-2x)^{-3} (2-3x)\\
        P_3(x) & = 3!(1-x)^{-4}(1-2x)^{-5} (5-14x+10x^2x)\\
        P_4(x) & = 4!(1-x)^{-5}(1-2x)^{-7} (14-56x+76x^2-35x^3)
      \end{split}
    \end{equation}
    and in general 
    \begin{equation}
      \label{eq:pkvalgen}
      P_k(x) = k!(1-x)^{-1-k}(1-2x)^{1-2k} Q_k(x)
    \end{equation}
    where \(Q_k(x)\) is a primitive polynomial of degree \(k-1\), with
    \(Q_k(1) = (-1)^{k+1}\).

    \begin{theorem}\label{thm:h2final}
      Put
      \begin{equation}
        \label{eq:PP}
        P(x,y) = \sum_{i,j \ge 0} c_{i,j}^{(2)} x^i \frac{y^j}{j!}
      \end{equation}
Then
\begin{equation}
  \label{eq:Pis}
  P(x,y) = \frac{2}{1+\sqrt{1 - 4(y + x - x^2)}}
\end{equation}
    \end{theorem}
    \begin{proof}
      We get from the recurrence relation \eqref{eq:Pk} that
      \begin{equation}
        \label{eq:diffeq1}
        (1-2x) \frac{\partial P}{\partial y} = \frac{\partial P}{\partial x}
      \end{equation}
      Furthermore, \(P_0(x)=P(x,0)=(1-x)^{-1}\). 
      The proposed \(P(x,y)\) satisfies \eqref{eq:diffeq1}
      and the initial condition, so it is the solution.
    \end{proof}

    \begin{theorem}\label{thm:c22}
      With the notations above, 
      \begin{equation}
        \label{eq:c2}
        \begin{split}
          c_{0,n}^{(2)} &= \frac{(2n)!}{(n+1)!} \\
          c_{1,n}^{(2)} &= c_{0,n+1}^{(2)} =\frac{(2(n+1))!}{(n+2)!} \\
          c_{2,n}^{(2)} &= \frac{1}{2}c_{0,n+2}^{(2)} - c_{0,n+1}^{(2)} =
          \frac{1}{16}{\frac {\left (2\,{n}^{2}+6\,n+3\right ){2}^{2\,n+6} 
              \Gamma (n+3/2)}{\left (n+3\right )\sqrt {\pi }\left (n+2\right )}}
        \end{split}
      \end{equation}
    Thus, the sequences \((c_{0,n}^{(2)})_{n=0}^\infty\) 
    and \((c_{1,n}^{(2)})_{n=0}^\infty\) are translations of 
the sequence
\htmladdnormallink{A001761}{http://www.research.att.com/projects/OEIS?Anum=001761}
in The On-Line Encyclopedia of Integer Sequences \cite{Sloane} (OEIS).
    \end{theorem}
    \begin{proof}
      We have that 
      \begin{equation}
        \label{eq:Py}
        P(0,y) =  \frac{2}{1+\sqrt{1 - 4y }}
      \end{equation}
      which is the well-known ordinary generating function for the
      Catalan numbers. This proves the formula for \(c_{0,n}^{(2)}\).
      The recurrence \eqref{eq:hrec} gives \(c_{1,n}^{(2)} =
      c_{0,n+1}^{(2)}\) and \(c_{2,n}^{(2)} = \frac{1}{2} c_{0,n+2}^{(2)}
      - c_{0,n+1}^{(2)}\). Combining these two results, and
      simplifying, yields the theorem.
    \end{proof}
  \end{section}

  \begin{section}{Enumeration of unrestricted standard paths}
    \begin{definition}
    Let \(\NC_{n,i,j}\) denote the poset of compositions of
    \(n\)  having \(i\) parts of size 1 and
    \(j\) parts of size \(\ge 2\). 
     Let \(\gamma_{n,i,j}\) be the
    number of standard paths with endpoint in \(\NC_{n,i,j}\).
    \end{definition}

    Clearly, 
    \begin{equation}
      \label{eq:kis1}
      \begin{split}
      \gamma_{n,i,0} &= 
      \begin{cases}
        1 & \text{ if } i=n,  \\
        0 & \text{ otherwise}
      \end{cases}
      \\
      \gamma_{0,i,j} &= 
      \begin{cases}
        1 & \text{ if } i=j=0  \\
        0 & \text{ otherwise}
      \end{cases}
      \\
      \gamma_{1,i,j} &= 
      \begin{cases}
        1 & \text{ if } i=1,j=0  \\
        0 & \text{ otherwise}
      \end{cases}
      \\
      \gamma_{n,0,j} &= 
      \begin{cases}
        1 & \text{ if } n>1,j=1 \text{ or if } n=j=0  \\
        0 & \text{ otherwise}
      \end{cases}
    \end{split}
    \end{equation}
    so \(F(u,0,x) = \exp(ux)\), \(F(u,v,0)=1\), 
    \(\left.\frac{\partial F(u,v,x)}{\partial x}\right\rvert_{x=0}
    =u\),
    \(F(0,v,x)=v (\exp(x)-x)\).

    We will derive a recurrence relation for \(\gamma_{n,i,j}\).
    Suppose that \(j > 0\).
    A tableau  with
    \(i\) parts of 
    size \(1\) and \(j\) parts of size \(\ge 2\),
    can be obtained
    \begin{enumerate}[-]
    \item from a tableau with \(i-1\) parts of size 1 and \(j\) parts
      of size \(\ge 2\), by adding a part of size 1 to the left,
    \item or from a tableau with \(i-1\) parts of size 1 and \(j\) parts
      of size \(\ge 2\), by adding a part of size 1 to the right,
    \item or from a tableau with \(i+1\) parts of size 1 and \(j-1\) parts
      of size \(\ge 2\), by adding a box to a part of size 1,
    \item or from a tableau with \(i\) parts of size 1 and \(j\) parts
      of size \(\ge 2\), by adding a box to a part of size \(\ge 2\).
    \end{enumerate}
    This gives the recurrence
    \begin{equation}
      \label{eq:grec}
      \gamma_{n,i,j} = 
      2\gamma_{n-1,i-1,j} + (i+1)\gamma_{n-1,i+1,j-1} + j\gamma_{n-1,i,j}
    \end{equation}
    valid for \(j >0\). If \(j=0\) then 
    \begin{displaymath}
      \gamma_{n,i,0} = \gamma_{n-1,i-1,0} = 
      \begin{cases}
        1 & n=i\\
        0 & \text{otherwise}
      \end{cases}
    \end{displaymath}

    \begin{definition}
    Define
    \begin{equation}
      \label{eq:H}
      F(u,v,x) = \sum_{n \ge 0} \Bigl( \sum_{i,j}
      \gamma_{n,i,j} u^i v^j \Bigl) \frac{x^n}{n!}
    \end{equation}
    \end{definition}
    
    The first few terms of the series \(F(u,v,x)\) are:
    \begin{multline}
1 +ux + \left ({u}^{2}+v\right )\frac{x^2}{2!}
+ \left ({u}^{3}+v+4\,uv\right )\frac{x^3}{3!} \\
+\left ({u}^{4}+v+6\,uv+11\,{u}^{2}v+4\,{v}^{2}
\right )\frac{x^4}{4!} \\
+ \left ({u}^{5}+v+8\,uv+23\,{u}^{2}v+26\,{u}^{3}v+14\,{v}^{2}+30\,u{v}^
{2}\right )\frac{x^5}{5!} + \dots
    \end{multline}
Setting \(u=v=1\) we obtain
\begin{multline}
  1+x+2\,\frac{x^2}{2!}+6\,\frac{x^3}{3!}+23\,\frac{x^4}{4!}+103\,\frac{x^5}{5!}+518\,\frac{x^6}{6!}+2868\,
\frac{x^7}{7!}+17263\,\frac{x^8}{8!}+111925\,\frac{x^9}{9!} + \dots
\end{multline}

\begin{theorem}\label{thm:F}
  Put \(F(u,v,x) = H(u,v,x) + F(u,0,x) = H(u,v,x) + \exp(ux)\).
  Then \(H(u,v,x)\) is the solution to 
  \begin{equation}
    \label{eq:Hpde}
    \begin{split}
      \frac{\partial H}{\partial x} &= 
      v\Bigl[ \frac{\partial H}{\partial v} +
      \frac{\partial H}{\partial u} +
      x \exp(ux) \Bigr]
      + 2uH \\
      H(u,0,x) & = \exp(ux)- \exp(ux) = 0 \\
      H(u,v,0) & = 1 - \exp(0) = 0 \\
      \left.\frac{\partial H(u,v,x)}{\partial x}\right\rvert_{x=0} &=
      u - \left.\frac{\partial \exp(ux)}{\partial x}\right\rvert_{x=0}
      =0 \\
      H(0,v,x) &= v (\exp(x)-x)  
    \end{split}
  \end{equation}
\end{theorem}
\begin{proof}
  This follows at once from \eqref{eq:grec} and \eqref{eq:kis1}.
\end{proof}
  \end{section}

\bibliographystyle{plain}
\raggedright
\bibliography{journals,articles,snellman}

\begin{thebibliography}{1}

\bibitem{StPa}
Fran{\c{c}}ois Bergeron, Mireille Bousquet-M{\'e}lou, and Serge Dulucq.
\newblock Standard paths in the composition poset.
\newblock {\em Ann. Sci. Math. Qu\'ebec}, 19(2):139--151, 1995.

\bibitem{Sloane}
Neil J.~A. Sloane.
\newblock The on-line encyclopedia of integer sequences.
\newblock
  {\texttt{http://www.research.att.com/{\(\sim\)}njas/sequences/index.html}}.

\bibitem{Snellman:Ncterm}
Jan Snellman.
\newblock A poset classifying non-commutative term orders.
\newblock In {\em Discrete models: {C}ombinatorics, {C}omputation, and
  {G}eometry}, Discrete {M}athematics and {T}heorethical {C}omputer {S}cience
  {P}roceedings {\textbf{AA (DM-CCG)}}, pages 301--314, 2001.

\end{thebibliography}

\end{document}